\documentclass[12pt]{amsart}
\usepackage{amssymb}

\newtheorem{theorem}{Theorem}
\theoremstyle{plain}
\newtheorem{definition}{Definition}
\newtheorem{lemma}{Lemma}
\newtheorem{proposition}{Proposition}
\newtheorem{corollary}{Corollary}
\newtheorem{remark}{Remark}
\newtheorem{example}{Example}
\numberwithin{equation}{section}

\begin{document}
\title[Gevrey microlocal analysis]{A Gevrey microlocal analysis of multi-anisotropic differential operators}
\author{Chikh BOUZAR and Rachid CHAILI }
\address{Department of Mathematics, Oran-Essenia University, Algeria}
\email{bouzar@univ-oran.dz; bouzar@yahoo.com}
\address{Department of Mathematics, University of Sciences and Technology of Oran,
Algeria} \email{chaili@univ-usto.dz} 
\begin{abstract}
We give a microlocal version of the theorem of iterates in multi-anisotropic
Gevrey classes for multi-anisotropic hypoelliptic differential operators.
\end{abstract}

\maketitle

\section{Introduction}

A fundamental result of Gevrey microlocal regularity due to H\"{o}rmander
\cite{H} is
\begin{equation}
WF_{s}(u)\subset WF_{s}(P(x,D)u)\cup Char(P)\text{ ,}  \label{1.1}
\end{equation}
where $P(x,D)$ denotes a differential operator\ with analytic coefficients
in $\Omega $, $Char(P)$ its set of characteristic points $(x,\xi )\in \Omega
\times \mathbb{R}^{n}$ and $WF_{s}(u)$ is the Gevrey wave front of the
distribution $u\in \frak{D}^{\prime }\left( \Omega \right) .$

Let $WF_{s}(u,P(x,D)),$ see \cite{BC}, be the Gevrey wave front of the
distribution $u\in \frak{D}^{\prime }\left( \Omega \right) $ with respect to
the iterates of the operator $P(x,D),$ then the result (\ref{1.1}) is made
more precise by the following inclusion
\begin{equation}
WF_{s}(u)\subset WF_{s}(u,P(x,D))\cup Char(P)\text{ ,}  \label{1.2}
\end{equation}
since
\begin{equation}
WF_{s}(u,P(x,D))\subset WF_{s}(P(x,D)u)\text{ .}  \label{1.3}
\end{equation}

Various extensions and generalizations of results (\ref{1.1}) and (\ref{1.2}%
) have been obtained, according as one considers the classes of elliptic or
hypoelliptic differential operators or one considers different notions of
homogeneity associated to these classes of operators, see e. g. \cite{BC},
\cite{BCR}, \cite{C}, \cite{R}, \cite{Z} and \cite{Z2}.

In \cite{LR} a microlocal analysis of the so called inhomogeneous Gevrey
classes \cite{R}, see also \cite{CMR}, has been introduced. The $\varphi -$%
inhomogeneous Gevrey wave front of a\ distribution $u\in \frak{D}^{\prime
}\left( \Omega \right) ,$ denoted $WF_{\varphi }(u),$ is defined with
respect to a weight function $\varphi .$

The method of Newton's polyhedron, see \cite{GV} or \cite{BBR}, permits to
approach differential operators with respect to their
multi-quasihomogeneity. In this situation the $\varphi -$inhomogeneous
Gevrey wave front $WF_{\varphi }(u)$ is characterized by a weight $\varphi $
equals to the function$\ \left| \xi \right| _{\mathbb{P}}$ defined by the
Newton's polyhedron $\mathbb{P}\ $of the operator $P(x,D)$ and it is denoted
by $WF_{s,\mathbb{P}}\left( u\right) .$

An interpretation of the $\varphi -$inhomogeneous Gevrey microlocal
analysis\ to the multi-anisotropic case is given in the paper \cite{C},
where a theorem in the spirit of the result (\ref{1.1}) for a class of
multi-quasihomogeneous hypoelliptic differential operators is obtained.

The aim of this paper is to obtain a result in the spirit of (\ref{1.2}) for
a class of multi-anisotropic hypoelliptic differential operators including
the classes of operators studied in \cite{BC}, \ \cite{BCh}, \cite{CH}, \cite
{C}, \cite{Ha}, \cite{Z} and \cite{Z2}. The section 2 is an adapted
modification of the $\varphi -$inhomogeneous Gevrey wave front of
Liess-Rodino, see \cite{LR} and \cite{C}, to our multi-anisotropic case in
the spirit of \cite{H} and \cite{Z}. In section 3 we introduce and study the
multi-anisotropic Gevrey wave front with respect to the iterates of an
operator $P(x,D)$ and its Newton's polyhedron $\mathbb{P},$\ denoted $WF_{s,%
\mathbb{P}}\left( u,P(x,D)\right) ,$ the following section 4 gives the
microlocal result of type (\ref{1.2}) for the studied class of\ differential
operators. This class is microlocally characterized by the following
definition.

\begin{definition}
Let $x_{0}\in \Omega ,\xi _{0}\in \mathbb{R}^{n}\backslash \left\{ 0\right\}
$\ and $P(x,D)$ be a\ differential operator with coefficients in the
anisotropic Gevrey class $G^{s,q}\left( \Omega \right) ,$ we say $(x_{0},\xi
_{0})\notin \sum_{\rho ,\delta ,s}^{\mu ,\mu ^{\prime },\mathbb{P}}(P)$ if
there exists an open neighbourhood $U\subset \Omega $ of $x_{0},$ an open $q$%
-quasiconic neighbourhood $\Gamma \subset \mathbb{R}^{n}\backslash \left\{
0\right\} $ of $\xi _{0}$\ and a constant $c>0$ such that $\forall \left(
x,\xi \right) \in U\times \Gamma ,$%
\begin{equation*}
\left\{
\begin{tabular}{l}
$\left| \xi \right| _{\mathbb{P}}^{\mu ^{\prime }}\leq c\left| P\left( x,\xi
\right) \right| ,$ \\
$\left| D_{x}^{\alpha }D_{\xi }^{\beta }P\left( x,\xi \right) \right| \leq
c^{\left| \alpha \right| +1}<\alpha ,q>^{s\mu <\alpha ,q>}\left| P\left(
x,\xi \right) \right| \left| \xi \right| _{\mathbb{P}}^{\delta \left| \alpha
\right| -\rho \left| \beta \right| },$%
\end{tabular}
\right.
\end{equation*}
where the numbers $\rho ,\delta ,\mu ^{\prime }$ and $\mu $ satisfy $0\leq
\delta <\rho \leq 1$\ and $\delta \mu <\mu ^{\prime }\leq \mu .\medskip $
\end{definition}

The principal result of this work is the following theorem.

\begin{theorem}
Let $u\in \frak{D}^{\prime }\left( \Omega \right) ,$ $P(x,D)$ a\
differential operator with coefficients in $G^{s,q}\left( \Omega \right) $
and $\rho ,\delta ,\mu ^{\prime },\mu $ such that $0\leq \delta <\rho \leq 1$%
\ and $\delta \mu <\mu ^{\prime }\leq \mu \medskip ,$ then\
\begin{equation}
WF_{s^{\prime },\mathbb{P}}\left( u\right) \subset WF_{s,\mathbb{P}}\left(
u,P\right) \text{ }\cup \sum\nolimits_{\rho ,\delta ,s}^{\mu ,\mu ^{\prime },%
\mathbb{P}}(P),  \label{1.5}
\end{equation}
where $s^{\prime }=\max \left( \frac{s\mu }{\mu ^{\prime }-\delta \mu },%
\frac{s}{\rho -\delta }\right) .$
\end{theorem}

\section{Multi-anisotropic Gevrey wave front}

This section is an adaptation with a slight modification of the
inhomogeneous Gevrey microlocal analysis introduced in \cite{LR}, see also
\cite{R} and \cite{C}, to the multi-anisotropic case.

Let $\Omega $ be an open subset of $\mathbb{R}^{n}$ and let $P\left(
x,D\right) $ be a linear partial differential operator with coefficients in $%
C^{\infty }\left( \Omega \right) ,$%
\begin{equation*}
P\left( x,D\right) =\sum\limits_{\alpha \in \Lambda }a_{\alpha }\left(
x\right) D^{\alpha }\;,
\end{equation*}
where $\Lambda $ is a finite subset of $\mathbb{Z}_{+}^{n}.$

\begin{definition}
Let $x_{0}\in \Omega ,$ the Newton's polyhedron of the\ operator $P\left(
x,D\right) $ at the point $x_{0}$, denoted $\mathbb{P}\left( x_{0}\right) ,$
is the convex hull of
\begin{equation*}
\left\{ 0\right\} \cup \left\{ \alpha \in \mathbb{Z}_{+}^{n},a_{\alpha
}\left( x_{0}\right) \neq 0\right\} \text{ }.\text{ }
\end{equation*}
\end{definition}

\begin{remark}
A Newton's polyhedron $\mathbb{P}$ is always characterized by
\begin{equation*}
\mathbb{P}=\underset{a\in \mathcal{A}}{\cap }\left\{ \alpha \in \mathbb{R}%
_{+}^{n},<\alpha ,a>\leq 1\right\} \;,
\end{equation*}
where $\mathcal{A}\left( \mathbb{P}\right) $ is a finite subset of $\mathbb{R%
}^{n}.$
\end{remark}

\begin{definition}
The Newton's polyhedron $\mathbb{P}$ is said to be regular if for any $%
a=(a_{1},...,a_{j},...,a_{n})\in \mathcal{A}$ we have $a_{j}>0,\forall
j=1,..,n.$
\end{definition}

\begin{definition}
The operator $P\left( x,D\right) $ is said regular if it satisfies the
following conditions :

\begin{enumerate}
\item  $\mathbb{P}\left( x_{0}\right) =\mathbb{P}$, $\forall x_{0}\in \Omega
$ .

\item  $\mathbb{P}$ is a regular polyhedron.
\end{enumerate}
\end{definition}

\begin{remark}
In this paper we consider only regular operators.
\end{remark}

Let $\mathbb{P}$ be a regular polyhedron, we set
\begin{eqnarray*}
\mathcal{V}\left( \mathbb{P}\right) &=&\left\{
s^{0}=0,s^{1},..,s^{m}\right\} \text{ the set of\ the vertices of }\mathbb{P}%
\text{ \ .} \\
\mu _{j} &=&\max a_{j}^{-1}\text{ , }a\in \mathcal{A}\;. \\
\;\mu &=&\max \mu _{j}\text{ \ .} \\
q &=&\left( \frac{\mu }{\mu _{1}},..,\frac{\mu }{\mu _{n}}\right) \text{ \ .}
\\
k\left( \alpha \right) &=&\inf \left\{ t>0,t^{-1}\alpha \in \mathbb{P}%
\right\} \underset{a\in \mathcal{A}}{=\max }<\alpha ,a>\text{ .} \\
\;\left| \xi \right| _{\mathbb{P}} &=&\left( \sum\limits_{i=1}^{m}\left( \xi
^{2s^{i}}\right) ^{1/\mu }\right) ^{1/2}\text{ .} \\
\left| \xi \right| _{q} &=&\left( \sum\limits_{j=1}^{n}\left( \xi
_{j}\right) ^{2/q_{j}}\right) ^{1/2}\text{\ .}
\end{eqnarray*}
$\;$

\begin{definition}
\label{def2.1}Let $s\geq 1$ and $\mathbb{P}$ \ be a regular polyhedron, we
denote $G^{s,\mathbb{P}}\left( \Omega \right) $ the space of functions $u\in
C^{\infty }\left( \Omega \right) $ such that $\forall K$ compact of $\Omega ,
$ $\exists C>0,\forall \alpha \in \mathbb{Z}_{+}^{n},$%
\begin{equation}
\sup_{K}\left| D^{\alpha }u\right| \leq C^{\left| \alpha \right| +1}k\left(
\alpha \right) ^{s\mu k\left( \alpha \right) }\text{ .}  \label{2.1}
\end{equation}
\end{definition}

\begin{example}
If the operator $P$ is $l$-quasi-elliptic of order $m$, with the weight $%
l=(l_{1},...,l_{n}),$ so its Newton's polyhedron $\mathbb{P}$ is the simplex
of vertices $\left\{ 0,\;m_{j}e_{j},\;j=1,..,n\right\} ,$ which is obviously
regular. In this case the set $\mathcal{A}$ coincides with the vector $%
\sum\limits_{j=1}^{n}m_{j}^{-1}e_{j},$ and we have $\mu _{j}=m_{j},\;\mu =m,$
$\;l=q=(\frac{m}{m_{1}},..,\frac{m}{m_{n}}).$ If $\alpha \in \mathbb{Z}%
_{+}^{n},$ then $k\left( \alpha \right) =m^{-1}<\alpha ,q>$ and we obtain $%
G^{s,\mathbb{P}}\left( \Omega \right) =G^{s,q}\left( \Omega \right) $ the
anisotropic Gevrey space, i.e. the space of functions $u\in C^{\infty
}\left( \Omega \right) $ such that $\forall K$ compact of $\Omega ,$ $%
\exists C>0,\forall \alpha \in \mathbb{Z}_{+}^{n},$%
\begin{equation*}
\sup_{K}\left| D^{\alpha }u\right| \leq C^{\left| \alpha \right| +1}\alpha
_{1}!^{q_{1}}...\alpha _{n}!^{q_{n}}\text{ \ }.
\end{equation*}
\end{example}

The following lemma, obtained in \cite{C}, gives the existence of a
truncation sequence, following the fundamental lemma 2.2 of \cite{H} in the
multi-anisotropic case. The quasihomogeneous case is a result of \cite[lemma
1.2]{Z}.

\begin{lemma}
Let $K$ be a compact set of $\mathbb{R}^{n}$ and let $s\geq 1,$ then there
exists a sequence $\left( \chi _{N}\right) \subset C_{0}^{\infty }\left(
\mathbb{R}^{n}\right) $ such that $\chi _{N}=1$ on $K$ and
\begin{equation}
\left| D^{\alpha }\chi _{N}\right| \leq C\left( CN^{s\mu }\right) ^{<\alpha
,a>}\text{ \ if}<\alpha ,a>\leq N,\forall a\in \mathcal{A},\;N=1,2,...\text{
\ .}  \label{2.2}
\end{equation}
\end{lemma}

A characterization of $G^{s,\mathbb{P}}\left( \Omega \right) $ using\ the
Fourier transform is given by the following theorem.

\begin{theorem}
Let $x_{0}\in \Omega $ and $u\in \frak{D}^{\prime }\left( \Omega \right) ,$
then $u$ is $G^{s,\mathbb{P}}$ in a neighbourhood of $x_{0}$ if, and only if
there exists a neighbourhood $U$\ of $x_{0}$ and a sequence $\left(
u_{N}\right) $ in $\mathcal{E}^{\prime }\left( \Omega \right) $ such that

i) $u_{N}=u$ \ in $U,\;N=1,2,...$ .

ii) $u_{N}$ is bounded in $\mathcal{E}^{\prime }\left( \Omega \right) $ .

iii) $\left| \widehat{u_{N}}\left( \xi \right) \right| \leq C\left( \frac{%
CN^{s}}{\left| \xi \right| _{\mathbb{P}}}\right) ^{\mu N},\;N=1,2,...$ .
\end{theorem}

\begin{proof}
See \cite{LR} and \cite{C}.
\end{proof}

We give now a\ microlocalization of the definition \ref{def2.1}. It is an
adapted modification , in the spirit of \cite{H} and \cite{Z}, of the $%
\varphi -$inhomogeneous Gevrey wave front of Liess-Rodino, see \cite{LR} and
\cite{C}, to our multi-anisotropic case. It coincides exactly with the
classical definition of the quasihomogeneous case.

\begin{definition}
\label{def2.2}Let $x_{0}\in \Omega ,$ $\xi _{0}\in \mathbb{R}^{n}\backslash
\left\{ 0\right\} $ and $u\in \frak{D}^{\prime }\left( \Omega \right) ,$ we
say that $u$ is $G^{s,\mathbb{P}}-$microregular at $\left( x_{0},\xi
_{0}\right) ,$ we denote $\left( x_{0},\xi _{0}\right) \notin WF_{s,\mathbb{P%
}}\left( u\right) ,$ if there exists $C>0,$ a neighbourhood $U$ of $x_{0}$
in $\Omega ,$ a $q$-quasiconic neighbourhood $\Gamma $ of $\xi _{0}$ in $%
\mathbb{R}^{n}\backslash \left\{ 0\right\} $ and a sequence $\left(
u_{N}\right) \subset \mathcal{E}^{\prime }\left( \Omega \right) $ such that

i) $u_{N}=u$ \ \ \ in $\ \ U,\;N=1,2,...$ .

ii) $u_{N}$ is bounded in $\mathcal{E}^{\prime }\left( \Omega \right) $ .

iii) $\left| \widehat{u_{N}}\left( \xi \right) \right| \leq C\left( \frac{%
CN^{s}}{\left| \xi \right| _{\mathbb{P}}}\right) ^{\mu N},\;N=1,2,...,\;\xi
\in \Gamma $ .
\end{definition}

We recall that a subset $\Gamma \subset \mathbb{R}^{n}$ is said $q$%
-quasiconic if
\begin{equation*}
\forall \xi \in \Gamma ,\forall t>0:\;\left( t^{q_{1}}\xi
_{1},...,t^{q_{n}}\xi _{n}\right) \in \Gamma .
\end{equation*}

\begin{remark}
The definition \ref{def2.2}\ coincides exactly with the quasihomogeneous
case, see \cite{Z}, if the polyhedron $\mathbb{P}$\ is the simplex of
vertices $\left\{ 0,\;m_{j}e_{j},\;j=1,..,n\right\} .$
\end{remark}

Using the truncation sequence $\left( \chi _{N}\right) $ we obtain the
following lemma, see \cite{C}.

\begin{lemma}
\label{lem2.2}Let $u\in \frak{D}^{\prime }\left( \Omega \right) $ and $%
\left( x_{0},\xi _{0}\right) \notin $ $WF_{s,\mathbb{P}}\left( u\right) $
and let $U,\Gamma $ be as in definition \ref{def2.2}. If $K$ is a compact
neighbourhood of $x_{0}$ in $U,$ $F$ is a $q$-quasiconic compact
neighbourhood of $\xi _{0}$ in $\Gamma $ and $\left( \chi _{N}\right)
\subset C_{0}^{\infty }\left( U\right) $ equal to $1$ on $K$ satisfying $%
\left( \ref{2.2}\right) ,$ then there exists $p_{0}\in \mathbb{Z}%
_{+},N_{0}\in \mathbb{Z}_{+}$ such that the sequence $\left( \chi
_{p_{0}N+N_{0}}u\right) $ satisfies i)-iii) in $K$ and $F$.
\end{lemma}

We define the $G^{s,\mathbb{P}}-$singsupp$\left( u\right) $ as the
complementary of the biggest open subset of $\Omega $ where $u$ is $G^{s,%
\mathbb{P}}.$ The relation between the multi-anisotropic Gevrey wave front
and the multi-anisotropic Gevrey singular support is given by the following
proposition.

\begin{proposition}
Let $u$ be a distribution in $\Omega $, then the projection of $WF_{s,%
\mathbb{P}}\left( u\right) $ on $\Omega $\ is the $G^{s,\mathbb{P}}-$singsupp%
$\left( u\right) .$
\end{proposition}

\begin{proof}
It follows the similar proof of \cite{C}.
\end{proof}

The microlocal property of the differential operator $P\left( x,D\right) $
with respect to the $G^{s,\mathbb{P}}-$wave front $WF_{s,\mathbb{P}}\left(
u\right) $\ is given by the following theorem.

\begin{theorem}
\label{th2.2}Let $u\in \frak{D}^{\prime }\left( \Omega \right) $ and $%
P\left( x,D\right) $ be a differential operator with coefficients in $%
G^{s,q}\left( \Omega \right) $, then
\begin{equation}
WF_{s,\mathbb{P}}\left( Pu\right) \subset WF_{s,\mathbb{P}}\left( u\right) .
\label{2.3}
\end{equation}
\end{theorem}

\begin{proof}
See \cite{LR} and \cite{R}.
\end{proof}

\begin{remark}
The product of two functions of the space $G^{s,\mathbb{P}}\left( \Omega
\right) $ does not belong in general to $G^{s,\mathbb{P}}\left( \Omega
\right) ,$ but if $g\in G^{s,q}\left( \Omega \right) $ and $f\in G^{s,%
\mathbb{P}}\left( \Omega \right) ,$ then $gf\in G^{s,\mathbb{P}}\left(
\Omega \right) ,$ see \cite{HM}. This justifies the optimal choice of the
regularity of the coefficients of the operator $P\left( x,D\right) .$
\end{remark}

\section{Multi-anisotropic Gevrey wave front with respect to the iterates of
a differential operator}

The Gevrey microlocal analysis with respect to the iterates of a
differential operator has been introduced for the first time by P. Bolley
and J. Camus in \cite{BC} in the homogeneous case. L. Zanghirati in \cite{Z}
has adapted it to the quasihomogeneous case. The aim of this section is to
extend this analysis to the multi-quasihomogeneous case.

\begin{definition}
\label{def3.1}Let $r\in \mathbb{R}$ and$\;s\geq 1,$ we denote $G_{r}^{s,%
\mathbb{P}}\left( \Omega ,P\right) $ the space of distributions $u\in \frak{D%
}^{\prime }\left( \Omega \right) $ such that $\forall K\;compact\;of\;\Omega
,\exists C>0,\forall N\in \mathbb{Z}_{+},$%
\begin{equation*}
\left\| P^{N}u\right\| _{H^{r}\left( K\right) }\leq C\left( CN^{s}\right)
^{\mu N}.
\end{equation*}
The space of Gevrey vectors of the operator $P$ is by definition
\begin{equation*}
G^{s,\mathbb{P}}\left( \Omega ,P\right) =\underset{r\in \mathbb{R}}{\cup }%
G_{r}^{s,\mathbb{P}}\left( \Omega ,P\right) .
\end{equation*}
\end{definition}

The space of Gevrey vectors $G^{s,\mathbb{P}}\left( \Omega ,P\right) $ of
the operator $P$ is described with the help of the Fourier transform in the
following lemma.

\begin{lemma}
\label{lem3.1}Let $x_{0}\in \Omega $ and $u\in \frak{D}^{\prime }\left(
\Omega \right) ,$ then $u\in $\ $G^{s,\mathbb{P}}\left( V,P\right) $ for a
neighbourhood $V$ of $x_{0}$ if, and only if, there exists a neighbourhood $U
$\ of $x_{0},U\subset V,$ $C>0,$ $M\in \mathbb{R}$ and a sequence $\left(
f_{N}\right) $ in $\mathcal{E}^{\prime }\left( V\right) $ such that

l) $f_{N}=P^{N}u$ \ in $U,\;N=0,1,...$ .

ll) $\left| \widehat{f_{N}}\left( \xi \right) \right| \leq C\left(
CN^{s}\right) ^{\mu N}\left( 1+\left| \xi \right| \right) ^{M},\;\xi \in
\mathbb{R}^{n},N=0,1,...$ .
\end{lemma}

\begin{proof}
It follows the proof of proposition 1.4 of \cite{BC}.
\end{proof}

The following technical lemma is important for the sequel.

\begin{lemma}
\label{lem3.2}Let $K$ be a compact subset of $\Omega $ and $\left( \chi
_{N}\right) $ a sequence in $C_{0}^{\infty }\left( \mathbb{R}^{n}\right) $
satisfying $\left( \ref{2.2}\right) $ and $\frak{F}$ a subset of $%
G^{s,q}\left( \Omega \right) $\ such that
\begin{equation*}
\exists C>0,\;\forall a\in \mathcal{A},\;\forall v\in \frak{F}%
,\;\;\sup_{K}\left| D^{\alpha }v\right| \leq C\left( C<\alpha ,a>^{s}\right)
^{\mu <\alpha ,a>},
\end{equation*}
then $\exists C_{1}>0,\;\forall v_{1},..,v_{j-1}\in \frak{F},\;\forall
\alpha ^{1},..,\alpha ^{j}\in \mathbb{Z}_{+}^{n},\forall a^{1},..,a^{j}\in
\mathcal{A},$\newline
$<\alpha ^{1},a^{1}>+..+<\alpha ^{j},a^{j}>\leq N,\ $we have
\begin{eqnarray*}
&&\;\;\;\;\;\;\;\;\;\;\;\;\sup_{K}\left| D^{\alpha ^{1}}v_{1}D^{\alpha
^{2}}v_{2}...D^{\alpha ^{j-1}}v_{j-1}D^{\alpha ^{j}}\chi _{N}\right|  \\
&\leq &C_{1}^{N+1}\left( \left( <\alpha ^{1},a^{1}>\right) ^{\mu <\alpha
^{1},a^{1}>}\left( <\alpha ^{2},a^{2}>\right) ^{\mu ^{2}<\alpha
^{2},a^{2}>}...\left( <\alpha ^{j},a^{j}>\right) ^{\mu ^{2}<\alpha
^{j},a^{j}>}\right) ^{s}.
\end{eqnarray*}
\end{lemma}

\begin{proof}
It is sufficient to see at first that $G^{s,q}\left( \Omega \right) \subset
G^{s,\mu a}\left( \Omega \right) ,\forall a\in \mathcal{A}$ and for any $%
a,b\in \mathcal{A},$ we have
\begin{equation*}
\left( \mu <\alpha ,a>\right) ^{\left( \mu <\alpha ,a>\right) }\leq \left(
\mu ^{2}<\alpha ,b>\right) ^{\left( \mu ^{2}<\alpha ,b>\right) },\;\alpha
\in \mathbb{Z}_{+}^{n},
\end{equation*}
and apply after lemma 2.3 of \cite{Z} or adapt lemma 5.3 of \cite{H}.
\end{proof}

Thanks to the truncation sequence $\left( \chi _{N}\right) ,$\ if $u\in
\frak{D}^{\prime }\left( \Omega \right) ,$ the sequence $u_{N}=\chi _{N}u$
is bounded in $\mathcal{E}^{\prime }\left( \Omega \right) $ and then $%
\exists C>0,$ $\left| \widehat{u_{N}}\left( \xi \right) \right| \leq C\left(
1+\left| \xi \right| \right) ^{M},\;\xi \in \mathbb{R}^{n},N\in \mathbb{Z}%
_{+}.$ In the problem of iterates this property is precised by the following
result.

\begin{lemma}
\label{lem3.3}Let $K$ be a compact subset\ of $\Omega $ and let $\left( \chi
_{N}\right) $ be a sequence in $C_{0}^{\infty }\left( K\right) $ satisfying $%
\left( \ref{2.2}\right) ,$ then $\forall u\in \frak{D}^{\prime }\left(
\Omega \right) ,\exists p_{0}>0,\forall p>p_{0},\forall r\in \mathbb{Z}_{+},$
the sequence $f_{N}=\chi _{pN+r}P^{N}$ satisfies
\begin{equation}
\widehat{f_{N}}\left( \xi \right) \leq C\left( C\left( N^{s\mu }+\left| \xi
\right| _{\mathbb{P}}\right) \right) ^{\mu N+M},\;\xi \in \mathbb{R}%
^{n},N\in \mathbb{Z}_{+},  \label{3.2}
\end{equation}
\end{lemma}

\begin{proof}
It does not differ substantially from its quasihomogeneous similar lemma 2.4
of \cite{Z}.
\end{proof}

The belonging to the space $G^{s,\mathbb{P}}\left( \Omega ,P\right) $ is
microlocally characterized by the following definition.

\begin{definition}
\label{def3.2}Let $u\in \frak{D}^{\prime }\left( \Omega \right) ,$ $%
(x_{0},\xi _{0})\in \Omega \times \mathbb{R}^{n}\backslash \left\{ 0\right\}
$ and $P(x,D)$ be a differential operator with coefficients in $%
G^{s,q}\left( \Omega \right) $. We say that $u$ is $G^{s,\mathbb{P}}-$%
microregular with respect to the iterates of $P(x,D)$ at $\left( x_{0},\xi
_{0}\right) $, we denote $\left( x_{0},\xi _{0}\right) \notin WF_{s,\mathbb{P%
}}\left( u,P\right) $, if there exists $C>0,M\in $ $\mathbb{R},$ a
neighbourhood $U$ of $x_{0}$ in $\Omega ,$ a $q$-quasiconic\ neighbourhood $%
\Gamma $ of $\xi _{0}$ in $\mathbb{R}^{n}\backslash \left\{ 0\right\} $ and
a sequence $\left( f_{N}\right) \subset \mathcal{E}^{\prime }\left( \Omega
\right) $ such that

j) $f_{N}=P^{N}u$ \ in $U,\;N\in \mathbb{Z}_{+}\;$\ .

jj)$\left| \widehat{f_{N}}\left( \xi \right) \right| \leq C\left( C\left(
N^{s\mu }+\left| \xi \right| _{\mathbb{P}}\right) \right) ^{\mu N+M},\;\xi
\in \mathbb{R}^{n},N\in \mathbb{Z}_{+}$ \ .

jjj)$\left| \widehat{f_{N}}\left( \xi \right) \right| \leq C\left(
CN^{s}\right) ^{\mu N}\left( 1+\left| \xi \right| \right) ^{M},\;\xi \in
\Gamma ,N\in \mathbb{Z}_{+}$ \ .
\end{definition}

The following proposition gives the link between the $G^{s,\mathbb{P}}-$%
singularities of a distribution $u\in \frak{D}^{\prime }\left( \Omega
\right) $ with respect to the iterates of $P(x,D)$\ and the wave front $%
WF_{s,\mathbb{P}}\left( u,P\right) .$

\begin{proposition}
Let $u\in \frak{D}^{\prime }\left( \Omega \right) $ and $P(x,D)$ be a
differential operator with coefficients in $G^{s,q}\left( \Omega \right) $,
then the projection of $WF_{s,\mathbb{P}}\left( u,P\right) $ on $\Omega $ is
the complementary of the biggest open subset $\Omega ^{\prime }$ of $\Omega $%
\ where $u\in G^{s,\mathbb{P}}\left( \Omega ^{\prime },P\right) .$
\end{proposition}

\begin{proof}
It follows the steps of the proofs of the classical theorems in the
homogeneous case, see \cite{BC}, and the quasihomogeneous case, see \cite{Z}%
, and makes use essentially of the following lemma.

\begin{lemma}
\label{lem3.4}Let $u\in \frak{D}^{\prime }\left( \Omega \right) $ and $%
\left( x_{0},\xi _{0}\right) \notin WF_{s,\mathbb{P}}\left( u,P\right) ,$ $U$
and $\Gamma $ be as in the definition \ref{def3.2}, $K$ a compact
neighbourhood of $x_{0}$ in $U,$ $F$ be a $q$-quasiconic compact
neighbourhood of $\xi _{0}$ in $\Gamma $ and $\left( \chi _{N}\right)
\subset C_{0}^{\infty }\left( U\right) $ be a sequence equals to $1$ on $K$
satisfying $\left( \ref{2.2}\right) ,$ then there exists $p_{0}\in \mathbb{Z}%
_{+},N_{0}\in \mathbb{Z}_{+}$ such that the sequence $\left( \chi
_{p_{0}N+N_{0}}P^{N}u\right) $ satisfies jjj) in $F$.
\end{lemma}
\end{proof}

The microlocal property of the operator $P(x,D)$\ with respect to the wave
front $WF_{s,\mathbb{P}}\left( u,P\right) $\ is the following result.

\begin{theorem}
\label{th3.1}Let $u\in \frak{D}^{\prime }\left( \Omega \right) $ and $P(x,D)$
be a differential operator with coefficients in $G^{s,q}\left( \Omega
\right) $, then
\begin{equation}
WF_{s,\mathbb{P}}\left( u,P\right) \subset WF_{s,\mathbb{P}}\left( Pu\right)
\subset WF_{s,\mathbb{P}}\left( u\right)
\end{equation}
\end{theorem}

\begin{proof}
Suppose that $\left( x_{0},\xi _{0}\right) \notin WF_{s,\mathbb{P}}\left(
u\right) $, then there exists a neighbourhood $U$ of $x_{0},$ a $q$%
-quasiconic neighbourhood $\Gamma $ of $\xi _{0}$ and a bounded sequence $%
\left( u_{N}\right) $ in $\mathcal{E}^{\prime }\left( \Omega \right) $ such
that $u_{N}=u$ in $U$ and $\left| \widehat{u_{N}}\left( \xi \right) \right|
\leq C\left( \frac{CN^{s}}{\left| \xi \right| _{\mathbb{P}}}\right) ^{\mu
N},\;N=1,2,..,\;\xi \in \Gamma .$ Let $K$ be a compact neighbourhood of $%
x_{0}$ in $U,$ $F$ be a $q$-quasiconic compact neighbourhood of $\xi _{0}$
in $\Gamma $ and let $\left( \chi _{N}\right) \subset C_{0}^{\infty }\left(
U\right) $ equal to $1$ on $K$ satisfying $\left( \ref{2.2}\right) .$ Choose$%
\;p\geq p_{0}+N_{0}$ and set $f_{N}=\chi _{pN}P^{N}u,$ we will show that
this sequence satisfies jjj) since j) is true and jj) is fulfilled according
to lemma \ref{lem3.3}.\newline
We have
\begin{equation}
\widehat{f}_{N}\left( \xi \right) =\int e^{-i<x,\xi >}\chi
_{pN}P^{N}udx=\int u\;^{t}P^{N}\left( e^{-i<x,\xi >}\chi _{pN}\right) dx.
\label{3.6}
\end{equation}
Set $^{t}P\left( x,D\right) =\sum\limits_{\alpha \in \mathbb{Z}_{+}^{n}\cap
\mathbb{P}}a_{\alpha }^{\prime }\left( x\right) D^{\alpha }$ and let $%
0=k_{0}<k_{1}<..<k_{r}=1,$ be the elements of the set $\left\{ k=k\left(
\alpha \right) ,\;\alpha \in \mathbb{Z}_{+}^{n}\cap \mathbb{P}\right\} .$\
Then
\begin{equation*}
^{t}P\left( e^{-i<x,\xi >}\chi _{pN+r}\right) =e^{-i<x,\xi >}\left| \xi
\right| _{\mathbb{P}}^{\mu }R\chi _{pN+r}\;,
\end{equation*}
where $R\left( x,\xi ,D\right) =R_{0}+...+R_{r}$ and
\begin{equation*}
R_{l}\left( x,\xi ,D\right) =\sum\limits_{\alpha \in \mathbb{Z}_{+}^{n}\cap
\mathbb{P}}\sum\limits_{\substack{ \beta \leq \alpha  \\ k\left( \beta
\right) =k_{l}}}\left( -1\right) ^{\left| \beta \right| }a_{\alpha }^{\prime
}\left( x\right) \frac{\xi ^{\beta }}{\left| \xi \right| _{\mathbb{P}}^{\mu }%
}D^{\alpha -\beta }\;.
\end{equation*}
By iteration we find
\begin{equation}
^{t}P^{N}\left( e^{-i<x,\xi >}\chi _{pN}\right) =e^{-i<x,\xi >}\left| \xi
\right| _{\mathbb{P}}^{\mu N}R^{N}\chi _{pN}=e^{-i<x,\xi >}\left| \xi
\right| _{\mathbb{P}}^{\mu N}\sum\limits_{\substack{ 0\leq l_{i}\leq r \\ %
1\leq i\leq N}}R_{l_{1}}...R_{l_{N}}\chi _{pN}\;.  \label{3.7}
\end{equation}
\newline
Since the coefficients of $R_{l}$\ are in $G^{s,q}\left( \Omega \right) ,$ $%
\forall \xi \in \mathbb{R}^{n},$ then from lemma \ref{lem3.2}, we obtain for
$<\alpha ,a>\leq N,$ $a\in \mathcal{A}$,
\begin{equation*}
\left| D^{\alpha }R_{l_{1}}...R_{l_{N}}\chi _{pN+r}\right| \leq
C_{1}^{N+1}N^{s\left( \mu <\alpha ,a>+\mu ^{2}N-\sum\limits_{1\leq i\leq
N}\mu ^{2}k_{l_{i}}\right) }\left( \left| \xi \right| _{\mathbb{P}}^{\mu
}\right) ^{\sum\limits_{1\leq i\leq N}k_{l_{i}}-N},
\end{equation*}
since $\left| \xi ^{\beta }\right| \leq \left| \xi \right| _{\mathbb{P}%
}^{\mu k\left( \beta \right) },\;\forall \beta \in \mathbb{Z}_{+}^{n}.$ Then
for $\left| \xi \right| _{\mathbb{P}}\geq N^{s\mu },<\alpha ,a>\leq N,$ $%
a\in \mathcal{A}$, we get
\begin{equation}
\left| D^{\alpha }\left( R^{N}\chi _{pN}\right) \right| \leq
C_{2}^{N+1}N^{s\mu <\alpha ,a>}.  \label{3.8}
\end{equation}

From $\left( \ref{3.6}\right) ,\left( \ref{3.7}\right) ,\left( \ref{3.8}%
\right) $ and lemma \ref{lem2.2}, we obtain
\begin{equation*}
\left| \widehat{f}_{N}\left( \xi \right) \right| =\left| \left( \left(
R^{N}\chi _{pN}\right) u\right) \symbol{94}\left( \xi \right) \right| \leq
C\left( CN^{s}\right) ^{\mu N}\;,\xi \in F,\left| \xi \right| _{\mathbb{P}%
}\geq N^{s\mu },
\end{equation*}
so $\left( x_{0},\xi _{0}\right) \notin WF_{s,\mathbb{P}}\left( u,P\right) ,$
hence $WF_{s,\mathbb{P}}\left( u,P\right) \subset WF_{s,\mathbb{P}}\left(
u\right) .$

Since
\begin{equation*}
WF_{s,\mathbb{P}}\left( u,P\right) =WF_{s,\mathbb{P}}\left( Pu,P\right)
\subset WF_{s,\mathbb{P}}\left( Pu\right)
\end{equation*}
and $WF_{s,\mathbb{P}}\left( Pu\right) \subset WF_{s,\mathbb{P}}\left(
u\right) $, according to theorem \ref{th2.2}, so the proof of theorem \ref
{th3.1} is complete.
\end{proof}

\section{The multi-anisotropic Gevrey microlocal regularity}

We obtain in this section a result of Gevrey microlocal regularity for a
class of multi-anisotropic hypoelliptic differential operators characterized
by the following definition.

\begin{definition}
\label{def4.1}Let $x_{0}\in \Omega ,\xi _{0}\in \mathbb{R}^{n}\backslash
\left\{ 0\right\} $\ and $P(x,D)$ be a\ differential operator with
coefficients in $G^{s,q}\left( \Omega \right) ,$ we denote $(x_{0},\xi
_{0})\notin \sum_{\rho ,\delta ,s}^{\mu ,\mu ^{\prime },\mathbb{P}}(P)$ if
there exists an open neighbourhood $U\subset \Omega $ of $x_{0},$ an open $q$%
-quasiconic neighbourhood $\Gamma \subset \mathbb{R}^{n}\backslash \left\{
0\right\} $ of $\xi _{0}$\ and a constant $c>0$ such that $\forall \left(
x,\xi \right) \in U\times \Gamma ,$%
\begin{equation}
\left\{
\begin{tabular}{l}
$\left| \xi \right| _{\mathbb{P}}^{\mu ^{\prime }}\leq c\left| P\left( x,\xi
\right) \right| ,$ \\
$\left| D_{x}^{\alpha }D_{\xi }^{\beta }P\left( x,\xi \right) \right| \leq
c^{\left| \alpha \right| +1}<\alpha ,q>^{s\mu <\alpha ,q>}\left| P\left(
x,\xi \right) \right| \left| \xi \right| _{\mathbb{P}}^{\delta \left| \alpha
\right| -\rho \left| \beta \right| },$%
\end{tabular}
\right.   \label{4.1}
\end{equation}
where the numbers $\rho ,\delta ,\mu ^{\prime }$ and $\mu $ satisfy $0\leq
\delta <\rho \leq 1$\ and $\delta \mu <\mu ^{\prime }\leq \mu .$
\end{definition}

We need the following lemma which is a modification of the similar result of
\cite[lemme 3.8]{BC}.

\begin{lemma}
\label{lem3.8}Under the notations of definition \ref{def4.1}, if $\chi
_{N}\in C_{0}^{\infty }\left( U\right) $ satisfies $\left( \ref{2.2}\right) ,
$ so there exists $C>0$ such that for $\left( x,\xi \right) \in U\times
\Gamma ,h_{1},..,h_{j}\in \mathbb{Z}_{+},a\in \mathcal{A},<\alpha
^{1}+..+\alpha ^{j},a>\leq N,$ $\beta ^{1},..,\beta ^{j-1}\in \mathbb{Z}%
_{+}^{n}:$%
\begin{eqnarray*}
&&\left| D^{\alpha ^{1}}P^{h_{1}}P^{\left( \beta ^{1}\right) }...D^{\alpha
^{j-1}}P^{h_{j-1}}P^{\left( \beta ^{j-1}\right) }D^{\alpha
^{j}}P^{h_{j}}\chi _{N}\right|  \\
&\leq &C^{N+1+\left| h_{1}\right| +..+\left| h_{j}\right| }<\alpha ,a>^{s\mu
<\alpha ,a>}\left| P\left( x,\xi \right) \right| ^{h_{1}+..+h_{j}+j-1}\left|
\xi \right| _{\mathbb{P}}^{\delta \left| \alpha \right| -\rho \left| \beta
\right| },
\end{eqnarray*}
where $\alpha =\alpha ^{1}+..+\alpha ^{j},\;\beta =\beta ^{1}+..+\beta
^{j-1}.$
\end{lemma}

The principal result of this work is the following theorem.

\begin{theorem}
\label{princ}Let $\Omega $ be an open subset of $\mathbb{R}^{n}$, $u\in
\frak{D}^{\prime }\left( \Omega \right) $ and $P(x,D)$ be a differential
operator with coefficients in $G^{s,q}\left( \Omega \right) $ and let $\rho
,\delta ,\mu ^{\prime }$ and $\mu $ be real numbers satisfying $0\leq \delta
<\rho \leq 1$\ and $\delta \mu <\mu ^{\prime }\leq \mu ,$ then\
\begin{equation}
WF_{s^{\prime },\mathbb{P}}\left( u\right) \subset WF_{s,\mathbb{P}}\left(
u,P\right) \cup \sum\nolimits_{\rho ,\delta ,s}^{\mu ,\mu ^{\prime },\mathbb{%
P}}\left( P\right) ,  \label{4.2}
\end{equation}
where $s^{\prime }=\max \left( \frac{s\mu }{\mu ^{\prime }-\delta \mu },%
\frac{s}{\rho -\delta }\right) .$
\end{theorem}

\begin{proof}
Let $\left( x_{0},\xi _{0}\right) \notin $ $WF_{s,\mathbb{P}}\left(
u,P\right) \cup \sum\nolimits_{\rho ,\delta ,s}^{\mu ,\mu ^{\prime },\mathbb{%
P}}\left( P\right) ,$ then there exists $C>0,M\in $ $\mathbb{R},$ a
neighbourhood $U$ of $x_{0}$ in $\Omega ,$ a $q$-quasiconic\ neighbourhood $%
\Gamma $ of $\xi _{0}$ in $\mathbb{R}^{n}\backslash \left\{ 0\right\} $ and
a sequence $\left( f_{N}\right) \subset \mathcal{E}^{\prime }\left( \Omega
\right) $ such that the conditions j),jj) and jjj) of definition \ref{def3.2}
are fulfilled.\ Let $K$ be a compact neighbourhood of $x_{0}$ in $U,$ $F$ a $%
q-$quasiconic compact neighbourhood of $\xi _{0}$ in $\Gamma $ such that $%
\left( \ref{4.1}\right) \ $is hold and let $\chi _{N}\in C_{0}^{\infty
}\left( U\right) ,$\ $\chi _{N}=1$ on $K$ satisfying $\left( \ref{2.2}%
\right) $ and $p$ a large enough integer. Set $u_{N}=\chi _{pN}u$ and let's
prove that this sequence satisfies iii) since i) and ii) are fulfilled. We
write
\begin{equation*}
^{t}P\left( e^{-i<x,\xi >}w\right) =e^{-i<x,\xi >}\left( ^{t}P\left( x,-\xi
\right) \left( I-R\right) \right) w,
\end{equation*}
where
\begin{equation*}
-R\left( x,\xi ,D\right) =\sum\limits_{\beta \neq 0}\frac{1}{\beta !}\;\frac{%
^{t}P^{\left( \beta \right) }\left( x,-\xi \right) }{^{t}P\left( x,-\xi
\right) }D^{\beta }.
\end{equation*}
By iteration we get
\begin{equation*}
^{t}P^{N}\left( e^{-i<x,\xi >}w\right) =e^{-i<x,\xi >}\left( ^{t}P\left(
x,-\xi \right) \left( I-R\right) \right) ^{N}w.
\end{equation*}

The fact that we can divide by $^{t}P\left( x,-\xi \right) $ is due to the
following lemma which can easily be proved.

\begin{lemma}
If $(x_{0},\xi _{0})\notin \sum_{\rho ,\delta ,s}^{\mu ,\mu ^{\prime },%
\mathbb{P}}(P),$ then $(x_{0},-\xi _{0})\notin \sum_{\rho ,\delta ,s}^{\mu
,\mu ^{\prime },\mathbb{P}}(^{t}P).$
\end{lemma}

Set
\begin{equation*}
w_{N}=\sum\limits_{h_{1}+..+h_{N}\leq \left[ \frac{\mu ^{\prime }-\delta \mu
}{\rho -\delta }N\right] }R^{h_{1}}\left( ^{t}P\right)
^{-1}...R^{h_{N}}\left( ^{t}P\right) ^{-1}\chi _{pN},
\end{equation*}
where $^{t}P=\,^{t}P\left( x,-\xi \right) .$ Then this function satisfies
\begin{equation*}
\left( ^{t}P\left( I-R\right) \right) ^{N}w_{N}=\chi _{pN}-e_{N},
\end{equation*}
where
\begin{equation*}
e_{N}=\sum\limits_{j=1}^{N}\left( ^{t}P\left( I-R_{N}\right) \right)
^{N-j}\sum\limits_{h_{j}+..+h_{N}=\left[ \frac{\mu ^{\prime }-\delta \mu }{%
\rho -\delta }N\right] }\,^{t}PR^{h_{j}+1}\left( ^{t}P\right)
^{-1}R^{h_{j+1}}...R^{h_{N}}\left( ^{t}P\right) ^{-1}\chi _{pN}.
\end{equation*}
Hence
\begin{equation}
\widehat{u}_{N}\left( \xi \right) =\widehat{w_{N}f_{N}}\left( \xi \right) +%
\widehat{e_{N}u}\left( \xi \right) ,\;\xi \in F.  \label{4.3}
\end{equation}
We will estimate both terms of the second member of $\left( \ref{4.3}\right)
$. Let $a\in \mathcal{A\medskip }$ and $0=k_{0}<k_{1}<...<k_{r}=1,$ be the
elements of the set $\left\{ k=<\alpha ,a>,\alpha \in \mathbb{Z}_{+}^{n}\cap
\mathbb{P}\right\} .$ we write $R=R_{1}+...+R_{r}$ where
\begin{equation*}
-R_{l}\left( x,\xi ,D\right) =\sum\limits_{<\beta ,a>=k_{l}}\frac{1}{\beta !}%
\;\frac{^{t}P^{\left( \beta \right) }\left( x,-\xi \right) }{^{t}P\left(
x,-\xi \right) }D^{\beta },
\end{equation*}
then we have
\begin{equation*}
w_{N}=\sum\limits_{h_{1}+..+h_{N}\leq \left[ \frac{\mu ^{\prime }-\delta \mu
}{\rho -\delta }N\right] }\sum\limits_{\substack{ 1\leq 1_{j}\leq r \\ 1\leq
j\leq h_{1}}}...\sum\limits_{\substack{ 1\leq N_{j}\leq r \\ 1\leq j\leq
h_{N}}}\left( R_{1_{1}}...R_{1_{h_{1}}}\right) \left( ^{t}P\right)
^{-1}...\left( R_{N_{1}}...R_{N_{h_{N}}}\right) \left( ^{t}P\right)
^{-1}\chi _{pN}.
\end{equation*}

Since $<\alpha ,a>\leq \left| \alpha \right| \leq \mu <\alpha ,a>,$ so from
lemma \ref{lem3.8}, we have for $<\alpha ,a>$ $\leq N,a\in \mathcal{%
A\medskip },$%
\begin{eqnarray*}
&&\;\left| D^{\alpha }\left( R_{1_{1}}...R_{1_{h_{1}}}\right) \left(
^{t}P\right) ^{-1}...\left( R_{N_{1}}...R_{N_{h_{N}}}\right) \left(
^{t}P\right) ^{-1}\chi _{pN}\right|  \\
&\leq &\medskip \medskip C^{N+1}N^{s\mu \left( <\alpha
,a>+\sum\nolimits^{\ast }k_{l_{i}}\right) }\left| P\left( x,\xi \right)
\right| ^{-N}\left| \xi \right| _{\mathbb{P}}^{\left( \delta -\rho \right)
\sum^{\ast }k_{l_{i}}+\mu \delta <\alpha ,a>},\medskip \medskip
\end{eqnarray*}
where $\sum^{\ast }$ means the sum over $1\leq l\leq N,1\leq i\leq
h_{l},1\leq l_{i}\leq r.$ Since the number of terms in the sum $w_{N}$ is
bounded from above by $C_{0}^{N},$ so $\exists C>0$ such that, for $<\alpha
,a>\leq N,\xi \in F,\left| \xi \right| _{\mathbb{P}}^{\rho -\delta }\geq
N^{s\mu },$ we have
\begin{equation*}
\left| D^{\alpha }w_{N}\right| \leq C_{1}^{N+1}N^{s\mu <\alpha ,a>}\left|
\xi \right| _{\mathbb{P}}^{\left( \mu \delta -\mu ^{\prime }\right) N}\;,
\end{equation*}
from lemma \ref{lem3.4}, we obtain for $\left| \xi \right| _{\mathbb{P}%
}^{\rho -\delta }\geq N^{s\mu },$%
\begin{eqnarray}
\left| \widehat{w_{N}f_{N}}\left( \xi \right) \right|  &\leq &C_{1}\left(
C_{1}N^{s}\right) ^{\mu N}\left| \xi \right| ^{M}\left| \xi \right| _{%
\mathbb{P}}^{\left( \mu \delta -\mu ^{\prime }\right) N}  \notag \\
&\leq &C_{2}\left( \frac{C_{2}N^{\frac{s\mu }{\mu ^{\prime }-\delta \mu }}}{%
\left| \xi \right| _{\mathbb{P}}}\right) ^{\frac{\left( \mu \delta -\mu
^{\prime }\right) }{\mu }\mu N}\left| \xi \right| ^{M}.  \label{4.4}
\end{eqnarray}

By the same procedure in the estimate of $w_{N}$, we get for $e_{N},$
\begin{equation*}
\left| D^{\alpha }e_{N}\right| \leq C_{3}^{N+1}N^{s\mu <\alpha ,a>}\left(
\frac{N^{\frac{s}{\rho -\delta }}}{\left| \xi \right| _{\mathbb{P}}}\right)
^{\left( \rho -\delta \right) \mu N},\;<\alpha ,a>\leq N,\left| \xi \right|
_{\mathbb{P}}^{\rho -\delta }\geq N^{s\mu }.
\end{equation*}
\ Let $M_{1}$ be the order of the distribution $u$ in $K$, so
\begin{equation}
\left| \widehat{ue_{N}}\left( \xi \right) \right| \leq C_{3}^{\prime
N+1}\left| \xi \right| ^{M_{1}}\left( \frac{N^{\frac{s}{\rho -\delta }}}{%
\left| \xi \right| _{\mathbb{P}}}\right) ^{\left( \rho -\delta \right) \mu
N}.  \label{4.5}
\end{equation}
From $\left( \ref{4.3}\right) ,\left( \ref{4.4}\right) $ and $\left( \ref
{4.5}\right) $ we easily obtain that
\begin{equation*}
\left( x_{0},\xi _{0}\right) \notin WF_{s^{\prime },\mathbb{P}}\left(
u\right) .
\end{equation*}
\end{proof}

\section{Consequences}

This section gives some corollaries of the obtained result.

\begin{corollary}
If $P(x,D)$\ is a differential operator with analytic coefficients,
satisfying $\left( \ref{4.1}\right) $ with $\left| \xi \right| _{\mathbb{P}%
}=\left| \xi \right| ,$ then theorem \ref{princ} coincides with the
principal theorem 5.1 of Bolley-Camus \cite{BC}, i.e. $\forall s\geq 1,$
\begin{equation*}
WF_{s}\left( u\right) \subset WF_{s}\left( u,P\right) \cup
\sum\nolimits_{\rho ,\delta ,s}^{\mu ,\mu ^{\prime },\mathbb{P}}\left(
P\right) .
\end{equation*}
\end{corollary}

\begin{remark}
The results of \cite{CH} and \cite{Ha} can be included in this corollary.
\end{remark}

\begin{corollary}
If the differential operator $P(x,D)$\ is $q$-quasihomogeneous with
coefficients in $G^{s,q}\left( \Omega \right) ,$ then $\left| \xi \right| _{%
\mathbb{P}}=\left| \xi \right| _{q}$ and
\begin{equation*}
\sum\nolimits_{1,0,s}^{\mu ,\mu ,\mathbb{P}}\left( P\right) =\left\{ \left(
x,\xi \right) \in \Omega \times \mathbb{R}^{n}\backslash \left\{ 0\right\}
:P_{q}\left( x,\xi \right) =0\right\} ,
\end{equation*}
where $P_{q}\left( x,\xi \right) $ is the principal $q$- quasihomogeneous
part of $P\left( x,\xi \right) .$ Consequently theorem \ref{princ} coincides
with the principal theorem of \cite{Z}, i.e.\ $\forall s\geq 1$%
\begin{equation*}
WF_{s,q}\left( u\right) \subset WF_{s,q}\left( u,P\right) \cup \left\{
\left( x,\xi \right) \in \Omega \times \mathbb{R}^{n}\backslash \left\{
0\right\} :P_{q}\left( x,\xi \right) =0\right\} .
\end{equation*}
\end{corollary}

\begin{definition}
The operator $P(x,D)$ is said multi-quasielliptic in $\Omega ,$\ if it is
regular and $\forall x_{0}\in \Omega ,$%
\begin{equation*}
\exists C>0,\;\exists R\geq 0,\;\left( \left| \xi \right| _{\mathbb{P}%
}\right) ^{\mu \left( \mathbb{P}\right) }\leq C\left| P\left( x_{0},\xi
\right) \right| ,\;\forall \xi \in \mathbb{R}^{n},\left| \xi \right| \geq R.
\end{equation*}
\end{definition}

The multi-anisotropic Gevrey regularity of the solutions of
multi-quasielliptic differential equations, see \cite{Z2} and \cite{BCh}, is
obtained easily from the following microlocal result.

\begin{corollary}
Let $u\in \frak{D}^{\prime }\left( \Omega \right) $ and $P(x,D)$\ be a
multi-quasielliptic differential operator with coefficients in $%
G^{s,q}\left( \Omega \right) ,$ then $\forall s\geq 1,$ we have\
\begin{equation*}
WF_{s,\mathbb{P}}\left( u\right) =WF_{s,\mathbb{P}}\left( u,P\right) =WF_{s,%
\mathbb{P}}\left( Pu\right) .
\end{equation*}
\end{corollary}

\textbf{Acknowledgements:} {The authors thank Professor Luigi Rodino for the
useful discussions on the subject of this paper.}

\end{document}